C2. Cognitive Science (including Linguistics and Psychology)

**Quantification in ordinary language:
from a critique of set-theoretic approaches to a proof-theoretic proposal**

**Michele Abrusci**
abrusci@uniroma3.it
Università Roma Tre

**Christian Retoré**
christian.retore@labri.fr
Université de Bordeaux

We firstly show that the standard interpretation of natural quantification in mathematical logic does not provide a satisfying account of its original richness. In particular, it ignores the difference between generic and distributive readings. We claim that it is due to the use of a set theoretical framework. We therefore propose a proof theoretical treatment in terms of proofs and refutations. Thereafter we apply these ideas to quantifiers that are not first order definable like "*the majority of*".

# Quantification in ordinary language:
# from a critique of set-theoretic approaches to a proof-theoretic proposal

## Michele Abrusci, Christian Retoré

**Empirical data**

Quantifiers are quite common in ordinary language and even more in speciality languages like sciences, in particular mathematics... As corpus studies show, existential quantifiers, often introduced by a simple indefinite, are so common that they structure discourses as shown by *Discourse Representation Theory*. Corpora also show that universal quantifiers are rarer but they are not so uncommon, and finally the so-called generalised quantifiers are possibly even more frequent than the universal one, be they explicitly asserted or implicitly introduced by the definite article or some other determiner. We therefore disagree with the common opinion in natural language processing that quantifiers are rare and should be neglected.

**Usual quantifiers**

Usual quantifiers, as opposed to generalized quantifiers, are assumed to be properly described by usual syntax (proofs) and semantics (models). . The key results are completeness, compactness and Löwenheim-Skolem theorems. But one should be aware that it only works for first-order classical logic. When second or higher order is used, sub-Boolean algebras have to be considered. When other logic are considered even intuitionistic logic, much more complicated models have to used, e.g. sheaves of L-structures.

The completeness theorem for first order classical logic says that there is a coincidence between provability and truth, that is $\&_{x,I}P(x)$ in every domain I, coincides with having a proof of P(x) for a generic x. It should be observed that natural language makes a distinction between the two above statements: *each* and *all* rather concern the elements of a collection while *every*, *any* or bare plurals rather concern the generic element.[1] Quantifiers without restriction like *everyone* or *everything* are possible too: they apply to a single property. Existential quantifications are formulated with *some, there is, a* and without restriction by *someone, something.* The negation of existential propositions is expressed either by *no,* and without restriction by *no one, nothing, …* Observe that the negation of universal proposition (the fourth corner of Aristotle square of oppositions) is NOT lexicalised (in English and other languages). Psycholinguists know it is harder to grasp and often understood as a negative existential. It is simply formulated by Not All: *Not all laureates deserve their award.* One could say: *Some laureates do not deserve their award.* But in this second formulation, the focus, the theme is different.

---

[1] The collective reading, where the collection, is viewed as a single individual, is not addressed in this paper.

**Proving and refuting quantified statements**

On these very classical matters our point is that the standard viewpoint does not provide a correct account of ordinary quantifiers: all of them are interpreted as a logical formula and thereafter as a set inclusion, no matter whether they concern a generic element or each individual of a collection. We advocate for another viewpoint, which is *pragmatic* (on the philosophical side) and *proof theoretical* (on the mathematical side).
When can one assert a quantifier? How does one refute it? Proof theory, following Hilbert's ideas (tau and epsilon operators), or usual quantification rules follows the generic reading while the very same theory, using rules like the Omega rule, follows the distributive reading. The usual refutation of an universal statement, consists in finding a particular X without the property. In addition to this refutation, we also consider a conceptual refutation which appears in natural language but not in usual logic, which does not refer to any individual,
- (Reasoning about dogs…) Every dog may bite. (assertion)
- Rex doesn't. (individual refutation)
- Basset-hounds don't. (conceptual refutation)

**"The majority of" vs. "Most of"**

Tricky quantifiers, like *the majority of* or *most of* are not first-order definable. Articles and books of formal semantics assume these two quantifiers to be the same, or rather that *most of* should be interpreted as *the majority of.* But there is a linguistic evidence that they do not mean the same: *most of* combines the difficulty of more than $n\%$ with the problem of vague predicates (the $n$ is implicit and much larger than half).

So let us focus on *the majority of*, which is simpler. The expression
*($\phi$) the majority of M's* are A
is usually interpreted as $|A|>|M-A|$ where $|X|$ stands for the cardinal of $|X|$. This interpretation is refuted by a simple experiment, asking people about the following sentences:
1. The majority of natural numbers are prime.
2. The majority of natural numbers are not prime.

As far as cardinality is concerned there are as many prime numbers as non prime numbers! Hence, our understanding is not a question of cardinality. In textbooks, (2.) may be found, and the proof of it is that zero is the limit of the ratio between primes lower than n and n. Hence, the semantic interpretation of such quantifiers, should be more sophisticated than cardinality (possibly a notion of measure, which only applies to definable subsets, should be introduced). In a proof theoretical perspective, what could be rules for *more than half*? When do we assert *more than half*? How can we refute a *more than half* judgement?

We shall provide guidelines of a possible solution.
- To assert (ϕ), we think some knowledge is required ( e.g. we should know some basic properties $B_i$ which are true of the majority of the elements).
- More than half are A can be asserted whenever all the basic properties $B_i$ encounter A.
- Regarding refutation, there are two ways to refute such a statement:
    - either by saying that only less than half satisfy A (dual quantifier)
    - or by saying that a property P which is known to be true for the majority is incompatible with A (P does not meet A).

For having a generic reading, as opposed to the distributive one that we just gave, a notion of medium, average element, would be needed, or of probabilistic proof. For the time being, it is out of reach.

**References:**
Jean-Yves Girard *Le point aveugle.* Tomes I et II, Hermann, 2006.
Stanley Peters, Dag Westerstahl *Quantifiers in Language And Logic*, Oxford University Press, 2006.